\documentclass[12pt,reqno]{amsart}

\usepackage[letterpaper,hmargin=1.25in,vmargin=1.22in]{geometry}
\usepackage{amssymb}
\usepackage[all]{xy}
\usepackage[hyphens]{url}

\newtheorem{de}{Definition}[section]
\newtheorem{lem}[de]{Lemma}
\newtheorem{prop}[de]{Proposition}
\newtheorem{cor}[de]{Corollary}
\newtheorem{thm}[de]{Theorem}

\theoremstyle{remark}
\newtheorem{rem}[de]{Remark}
\newtheorem{ex}[de]{Example}

\def \N{\mathbb{N}}
\def \Z{\mathbb{Z}}

\def \R{\mathbb{R}}

\newcommand{\Diff}{\mathrm{Diff}}
\newcommand{\DVPB}{\mathrm{DVPB}}

\newcommand{\im}{\mathrm{Im}}
\newcommand{\dvs}{\mathrm{dvs}}
\newcommand{\Hom}{\mathrm{Hom}}

\newcommand{\cSD}{\mathcal{SD}}
\newcommand{\cSV}{\mathcal{SV}}

\newcommand{\Exterior}{\mathchoice{{\textstyle\bigwedge}}%
{{\bigwedge}}%
{{\textstyle\wedge}}%
{{\scriptstyle\wedge}}}

\title{Pushforward and smooth vector pseudo-bundles}

\author{Enxin Wu}
\email{exwu@stu.edu.cn}
\address{Department of Mathematics, Shantou University, Guangdong, P.R. China}

\date{\today}

\begin{document}

\begin{abstract}
 In this paper, we study a new operation named pushforward on diffeological vector pseudo-bundles, 
 which is left adjoint to the pullback. We show how to pushforward projective diffeological 
 vector pseudo-bundles to get projective diffeological vector spaces, producing many concrete new examples, 
 together with application to smooth splittings of some projective diffeological vector spaces related to geometry.
 This brings new objects to diffeology from classical vector bundle theory.
\end{abstract}

\makeatletter
\@namedef{subjclassname@2020}{%
  \textup{2020} Mathematics Subject Classification}
\makeatother


\maketitle

\section{Introduction}

Diffeological spaces are elegant generalizations of smooth manifolds, including infinite-dimensional spaces 
like mapping spaces and diffeomorphism groups, and singular spaces eg smooth manifolds with boundary or 
corners, orbifolds and irrational tori.

On diffeological spaces, one can still do some differential geometry and topology, such as differential forms 
and tangent bundles. These tangent bundles are in general no longer locally trivial. 
Instead, they are diffeological vector pseudo-bundles. We studied these objects and operations on them 
in~\cite{CWp}, on which the current paper is based.

On the other hand, the theory of diffeological vector spaces and their homological algebra is intimately 
related to analysis and geometry; see~\cite{W,CW16,CW21}. The projective objects there deserve special attention. 
However, in general neither is it easy to test whether a given diffeological vector space is projective or not, 
nor is it straightforward to construct many concrete projective objects. 

In this paper, we propose a way to use diffeological vector pseudo-bundles to study diffeological vector spaces.
We generalize some results of projective objects for diffeological vector spaces to such bundles. 
In particular, we show that every classical vector bundle is such a projective object.
We introduce a left adjoint called pushforward to the pullback on diffeological vector pseudo-bundles, 
and we show that the free diffeological vector space generated by a diffeological space has a canonical bundle-theoretical 
explanation, and that pushforward preserves projectives. In this way, we construct many concrete 
projective diffeological vector spaces from classical vector bundle theory, together with application of classical 
vector bundle theory to smooth splittings of some projective diffeological vector spaces.

Here is the structure of the paper. In Section~\ref{s:background}, we briefly review some necessary background. 
In Section~\ref{s:pushforward}, we introduce pushforward on diffeological vector pseudo-bundles.
Section~\ref{s:projective} contains three parts, including necessary and sufficient conditions of smooth splittings of 
short exact sequences of diffeological vector pseudo-bundles, examples and properties of the projective objects, 
and preservation of projectives by pushforward. In particular, we get many new examples of projective diffeological 
vector spaces from classical vector bundles. In Section~\ref{s:application}, we apply the established theory to 
smooth splittings of projective diffeological vector spaces. Readers interested in concrete examples are suggested to 
take a look at the last part of this section first.

\section{Background}\label{s:background}

We give a very brief review together with many related references in this section.

\begin{de}
 A \textbf{diffeological space} is a set $X$ together with a collection of maps $U \to X$ (called \textbf{plots})
 from open subsets $U$ of Euclidean spaces, such that
 \begin{enumerate}
  \item every constant map is a plot;
  \item The composite $V \to U \to X$ is a plot if the first map is smooth between open subsets of Euclidean spaces
        and the second one is a plot;
  \item $U \to X$ is a plot if there is an open cover of $U$ such that each restriction is a plot.
 \end{enumerate}

 A \textbf{smooth map} $X \to Y$ between diffeological spaces is a map which sends plots of $X$ to plots of $Y$.
 Diffeological spaces with smooth maps form a category denoted $\Diff$.
\end{de}

The idea of a diffeological space was introduced in~\cite{S}. \cite{I13} is currently the standard reference
for the subject. Also see~\cite[Section~2]{CSW} for a concise summary for the basics of diffeological spaces.

The category $\Diff$ has excellent properties. It contains the category of smooth manifolds as a full subcategory, 
and it is complete, cocomplete and cartesian closed. In particular, we have subspaces, quotient spaces 
and mapping spaces for diffeological spaces. Like charts for manifolds, we have various generating sets of plots for 
a diffeological space. Every diffeological space has a canonical topology called 
the $D$-topology; see~\cite{I85,CSW}. Every diffeological space has a tangent bundle; see~\cite{H,CW16,CW17}.
Diffeological vector spaces are the vector space objects in $\Diff$. 
Every vector space can be equipped with a smallest diffeology called the fine diffeology, 
making it a diffeological vector space; see~\cite{I07}. 
There are many other kinds of diffeological vector spaces in practice.
Hierachies of diffeological vector spaces were studied in~\cite{CW19}, 
and homological algebra of diffeological vector spaces, including free and projective objects, 
were introduced in~\cite{W}.

We recall the following concepts from~\cite{CWp}:

\begin{de}
A \textbf{diffeological vector pseudo-bundle} over a diffeological space $B$ is a smooth map $\pi:E \to B$
between diffeological spaces such that the following conditions hold:
\begin{enumerate}
 \item for each $b \in B$, $\pi^{-1}(b) =: E_b$ is a vector space;  
 \item the fibrewise addition $E \times_B E \to E$ and the fibrewise scalar multiplication 
  $\R \times E \to E$ are smooth;
 \item the zero section $\sigma:B \to E$ is smooth.
\end{enumerate}
\end{de}

\begin{de}
Given a diffeological space $B$, a \textbf{bundle map over $B$} is a commutative triangle
\[
 \xymatrix@C5pt{E_1 \ar[rr]^f \ar[dr]_{\pi_1} && E_2 \ar[dl]^{\pi_2} \\ & B,}
\]
where $\pi_1,\pi_2$ are diffeological vector pseudo-bundles over $B$, $f$ is smooth and for each $b \in B$, 
the restriction $f|_{E_{1,b}}:E_{1,b} \to E_{2,b}$ is linear. 

Such $f$ is called a \textbf{bundle subduction} (resp. \textbf{bundle induction}) \textbf{over $B$} 
if it is both a bundle map over $B$ and a subduction (resp. an induction), i.e., it is equivalent 
to a quotient map (resp. an inclusion of a subspace).
\end{de}

For a fixed diffeological space $B$, all diffeological vector pseudo-bundles over $B$ and bundle maps 
over $B$ form a category, denoted $\DVPB_B$. An isomorphism in $\DVPB_B$ is called a \textbf{bundle
isomorphism over $B$}. A bundle map over $B$ is a bundle isomorphism if and only if it is both a bundle induction and 
a bundle subduction over $B$.

\begin{de}
 A commutative square 
 \[
  \xymatrix{E \ar[r]^g \ar[d]_{\pi} & E' \ar[d]^{\pi'} \\ B \ar[r]_f & B'}
 \]
 in $\Diff$ with $\pi$ and $\pi'$ being diffeological vector pseudo-bundles, is called a \textbf{bundle map}, 
 if for each $b \in B$, $g|_{E_b}:E_b \to E'_{f(b)}$ is linear.

 A bundle map $(g,f)$ as above is called a \textbf{bundle subduction} if both $g$ and $f$ are subductions.
\end{de}

All diffeological vector pseudo-bundles and bundle maps form a category denoted $\DVPB$.

Note that diffeological vector pseudo-bundles are neither diffeological fibre bundles in~\cite{I85,I13}, 
nor diffeological fibrations in~\cite{CW14}. They were introduced to encode tangent bundles of diffeological 
spaces (\cite{CW16}).
Many operations on $\DVPB_B$ and $\DVPB$ were studied in~\cite{CWp}, such as direct product, direct sum, 
free diffeological vector pseudo-bundle induced by a smooth map, tensor product, and exterior product.

\section{Pushforward}\label{s:pushforward}

Recall from~\cite[Section~3.1]{CWp} that one can pullback diffeological vector pseudo-bundles via smooth maps, 
i.e., a smooth map $f:B \to B'$ induces a functor $f^*:\DVPB_{B'} \to \DVPB_B$ by pullback. 
Now we define a related operation as follows:

Given a smooth map $f:B \to B'$ and a diffeological vector pseudo-bundle $\pi:E \to B$, we define 
\begin{equation}\label{eq:pushforward}
 E' = \coprod_{b' \in B'} ( \bigoplus_{b \in f^{-1}(b')} E_b ).
\end{equation}
Note that when $f^{-1}(b') = \emptyset$, the term in the above bracket is $\R^0$. 
There are canonical maps $\pi_f:E' \to B'$ sending the fibre above $b'$ to $b'$, 
and $\alpha_f:E \to E'$ with $E_b \hookrightarrow \bigoplus_{\tilde{b} \in f^{-1}(f(b))} E_{\tilde{b}}$. 
We then have a natural commutative square
\[
 \xymatrix{E \ar[r]^{\alpha_f} \ar[d]_{\pi} & E' \ar[d]^{\pi_f} \\ B \ar[r]_f & B'.}
\]
Hence, we can equip $E'$ with the dvsification of the diffeology generated by the upper horizontal map $\alpha_f$
of the above square via~\cite[Proposition~3.3]{CWp}, making the right vertical map $\pi_f$ a diffeological vector pseudo-bundle over $B'$, and hence the above square becomes a bundle map from $\pi$ to $\pi_f$. 
(As a \textbf{warning}, each fibre of $E'$ may not be the direct sum of those of $E$ as diffeological vector 
spaces; see Proposition~\ref{prop:directsum}. Also notice that the notation $\alpha_f$ will be used later in the paper.)
More precisely, we have the following explicit description of a generating set of plots on $E'$:

\begin{lem}\label{lem:diffeology-pushforward}
 A plot on $E'$ is locally of one of the following forms:
 \begin{enumerate}
  \item $U \to E'$ defined by a finite sum $\sum_i \alpha_f \circ p_i$, where $p_i:U \to E$ are plots on $E$
        such that all $f \circ \pi \circ p_i$'s match;       
  \item the composite of a plot of $B'$ followed by the zero section $B' \to E'$.
 \end{enumerate}
\end{lem}
\begin{proof}
 This is straightforward from the description of dvsification in~\cite{CWp}.
\end{proof}

It is straightforward to check that we get a functor $f_*:\DVPB_B \to \DVPB_{B'}$, called the \textbf{pushforward 
of $f$}, and we write $E'$ above as $f_*(E)$.
Moreover, from the above lemma, we have 
\begin{enumerate}
 \item $f'_* \circ f_* = (f' \circ f)_*$ for any smooth maps $f:B \to B'$ and $f':B' \to B''$;
 \item $(1_B)_* =$ the identity on $\DVPB_B$.
\end{enumerate}

\begin{ex}
 Pushforward has been used implicitly in~\cite[Section~5]{CWp}. For example, $E_1$ and $E_2$
 in~\cite[Proposition~5.1]{CWp} are the pushforward of the tangent bundle $\R^2 \to \R$ along 
 the inclusions $\R \to X_g$ to the $x$-axis and the $y$-axis, respectively.
\end{ex}

Here is the key result for pushforward:

\begin{thm}\label{thm:adjoint}
 Given a smooth map $f:B \to B'$, we have an adjoint pair of functors
 \[
  f_*: \DVPB_B \rightleftharpoons \DVPB_{B'} : f^*.
 \]
\end{thm}
\begin{proof}
 We show that there is a natural bijection $\DVPB_B(E,f^*(E')) \cong \DVPB_{B'}(f_*(E),E')$. 
 Given a bundle map $E \to f^*(E')$ over $B$, we have $E_b \to E'_{f(b)}$ for each $b \in B$, 
 which induce $\bigoplus_{b \in f^{-1}(b')} E_b \to E'_{b'}$, and hence a map $f_*(E) \to E'$. 
 This is clearly a bundle map over $B'$. Conversely, given a bundle map $f_*(E) \to E'$ over $B'$, 
 we have a map $\bigoplus_{b \in f^{-1}(b')} E_b \to E'_{b'}$ for each $b' \in \im(f)$. 
 It then induces a map $E_b \to E'_{f(b)}$, which together give a map $E \to f^*(E')$. 
 It is straightforward to check that this is a bundle map over $B$. 
 These procedures are inverses to each other, and therefore we proved the desired result.
\end{proof}

We have the following bundle-theoretical explanation of a free diffeological vector space
introduced in~\cite{W}:

\begin{prop}\label{prop:free-geom}
 For any diffeological space $B$, the total space of the pushforward of the trivial bundle $B \times \R \to B$
 along the map $B \to \R^0$ is the free diffeological vector space $F(B)$.
\end{prop}
\begin{proof}
 This follows directly from the diffeology of the total space of the pushforward (see
 Lemma~\ref{lem:diffeology-pushforward}) and the 
 diffeology on free diffeological vector space (see proof of~\cite[Proposition~3.5]{W}).
\end{proof}

From~\cite[Section~3]{CWp}, we know that the usual operations on diffeological vector pseudo-bundles 
have the obvious diffeology on each fibre indicated by the operation. But pushforward is an exception, 
although it is expected so:

\begin{prop}\label{prop:directsum}
 Let $f:B \to B'$ be a smooth map, and let $E \to B$ be a diffeological vector pseudo-bundle.
 Then the diffeology on the fibre at $b'$ of the pushforward $f_*(E)$ has the direct sum diffeology of the
 diffeological vector spaces $E_b$'s with $f(b)=b'$ if and only if $f^{-1}(b')$ as a subspace of $B$ has
 the discrete diffeology.
\end{prop}
\begin{proof}
 This follows directly from Lemma~\ref{lem:diffeology-pushforward}.
\end{proof}

Here is the universal property for pushforward:

\begin{prop}\label{prop:universal}
 Given a bundle map 
 \[
  \xymatrix{E \ar[r]^f \ar[d]_{\pi} & E' \ar[d]^{\pi'} \\ B \ar[r]_g & B',}
 \]
 there exists a unique bundle map $\beta:g_*(E) \to E'$ over $B'$ such that $f = \beta \circ \alpha_g$.
\end{prop}
\begin{proof}
 This is clear by the construction of pushforward, or from the adjoint (Theorem~\ref{thm:adjoint}).
\end{proof}

Pushforward could send non-isomorphic bundles to isomorphic ones:

\begin{ex}\label{ex:cross}
 Write $B$ for the cross with the gluing diffeology, and write $B'$ for the cross with the subset diffeology 
 of $\R^2$. Then $B \to B'$ defined as the identity underlying set map is smooth, but its inverse is not;
 see~\cite[Example~3.19]{CW16}.
 We show below that the induced map $F(B) \to F(B')$ between the free diffeological vector spaces,
 which is identity for the underlying vector spaces,
 is indeed an isomorphism of diffeological vector spaces. 
 This means that the pushforward of the two trivial bundles $B \times \R \to B$ and $B' \times \R \to B'$ 
 along the maps $B \to \R^0$ and $B' \to \R^0$ are isomorphic, but clearly the two bundles are not.

 By definition of a free diffeological vector space, every plot $p:U \to F(B')$ can be locally written as 
 a finite sum $p(u) = \sum_i r_i(u) (p_{1i}(u),p_{2i}(u))$ for smooth maps $r_i,p_{1i},p_{2i}$ with codomain $\R$
 satisfying $p_{1i}(u) p_{2i}(u) = 0$ for all $u$. It is enough to show that $p$ can be viewed as a plot 
 of $F(B)$. This is the case since $(p_{1i}(u),p_{2i}(u))$ can be written as $(p_{1i}(u),0) + (0,p_{2i}(u)) - (0,0)$,
 each term viewed as a plot of $B$.
\end{ex} 

As a consequence of the above example, the canonical map $X \to F(X)$ from a diffeological space to the free 
diffeological vector space generated by it, is \emph{not} necessarily an induction.
Recall that $F(X)$ is the free vector space over $\R$ generated by the underlying set of $X$, 
and it is equipped with the smallest vector space diffeology such that the canonical map $X \to F(X)$ 
sending $x \in X$ to the base element $[x] \in F(X)$.
This observation is a bit surprising, isn't it?

On the other hand, we have

\begin{prop}
 The canonical map $X \to F(X)$ is an induction if and only if there exist a family of diffeological vector spaces 
 $\{V_i\}_{i \in I}$ such that the diffeology on $X$ is determined by the union of all $C^\infty(X,V_i)$, in the 
 sense that $U \to X$ is a plot if and only if the composite $U \to X \to V_i$ is smooth for every smooth map 
 $X \to V_i$. 
\end{prop}
In particular, for every Fr\"olicher space $X$ (i.e., the diffeology on $X$ is determined by $C^\infty(X,\R)$), 
the canonical map $X \to F(X)$ is an induction. This applies to $B'$ in Example~\ref{ex:cross}.
\begin{proof}
 This follows immediately from the universal property of the free diffeological vector space generated by a 
 diffeological space.
\end{proof}

\section{Projective diffeological vector pseudo-bundles}\label{s:projective}

\subsection{Enough projectives}

In this subsection, we will work in the category $\DVPB_B$ for a fixed diffeological space $B$. 
So we will omit the phrase `over $B$' in many places as long as no confusion shall occur.
Note that when we take $B = \R^0$, we recover the corresponding results for the category 
of diffeological vector spaces.

We first study smooth splittings of diffeological vector pseudo-bundles, 
which will be used later in the paper.

\begin{de}
 A diagram of morphisms
 \[
  \xymatrix{E_1 \ar[r]^f & E_2 \ar[r]^g & E_3}
 \]
 in $\DVPB_B$, 
 is called a \textbf{short exact sequence} if $f$ is a bundle induction, $g$ is a bundle subduction, and 
 \[
  \xymatrix{E_{1,b} \ar[r]^{f_b} & E_{2,b} \ar[r]^{g_b} & E_{3,b}}
 \]
 is exact (i.e., $\ker(g_b) = \im(f_b)$) for every $b \in B$.
\end{de}

As a direct consequence of the above definition, we have:

\begin{cor}
 Given a short exact sequence 
 \[
  \xymatrix{E_1 \ar[r] & E_2 \ar[r] & E_3}
 \]
 of diffeological vector pseudo-bundles over $B$, we have a bundle isomorphism $E_2/E_1 \cong E_3$ over $B$.
\end{cor}

The splitting of a short exact sequence goes as usual:

\begin{thm}\label{thm:splitses}
 Assume that 
 \[
  \xymatrix{E_1 \ar[r]^f & E_2 \ar[r]^g & E_3}
 \]
 is a short exact sequence of diffeological vector pseudo-bundles over $B$. 
 Then the following are equivalent:
 \begin{enumerate}
  \item there exists a bundle map $g':E_3 \to E_2$ over $B$ such that $g \circ g'= 1_{E_3}$;
  \item there exists a bundle map $f':E_2 \to E_1$ over $B$ such that $f' \circ f = 1_{E_1}$;
  \item there exists a bundle isomorphism $E_2 \to E_1 \oplus E_3$ over $B$
   making the following diagram commutative:
   \[
    \xymatrix{E_1 \ar[r]^f \ar[d]_{=} & E_2 \ar[r]^g \ar[d] & E_3 \ar[d]^{=} \\ 
              E_1 \ar[r]_-{i_1} & E_1 \oplus E_3 \ar[r]_-{p_2} & E_3}
   \]
 \end{enumerate}
\end{thm}
If any one of the conditions holds in the theorem, we say that the short exact sequence \textbf{splits smoothly}, 
and that $E_1$ (resp. $E_3$) is a \textbf{smooth direct summand} of $E_2$.
Although every short exact sequence of vector spaces splits, it is not the case in $\DVPB_B$, 
even when $B = \R^0$; see~\cite[Example~4.3]{W} or~\cite[Example~4.1]{CW19}.
\begin{proof}
 We show below that (1) $\Leftrightarrow$ (3), and (2) $\Leftrightarrow$ (3) can be proved similarly.

 (1) $\Rightarrow$ (3): since we have bundle maps $f:E_1 \to E_2$ and $g':E_3 \to E_2$, we define 
 $E_1 \oplus E_3 \to E_2$ by $(x_1,x_3) \mapsto f(x_1) + g'(x_3)$ for any $x_1 \in E_{1,b}$, $x_3 \in E_{3,b}$
 and $b \in B$. This is clearly a bundle map over $B$. Its inverse is given by 
 $x \mapsto (f^{-1}(x - g' \circ g(x)),g(x))$. It is straightforward to check that this is well-defined, 
 and it is smooth since $f$ is an induction.

 (3) $\Rightarrow$ (1): $g'$ is defined by the composite $\xymatrix{E_3 \ar[r]^-{i_2} & E_1 \oplus E_3 \ar[r]^-{\cong} & E_2}$.
 The rest are straightforward to check.
\end{proof}

Now we can define projective diffeological vector pseudo-bundles, and show that there are \emph{enough} 
such objects.

\begin{de}\label{de:projective}
 A diffeological vector pseudo-bundle $E \to B$ is called \textbf{projective} if for any bundle 
 subduction $f:E_1 \to E_2$ over $B$ and any bundle map $g:E \to E_2$ over $B$, there exists a bundle map
 $h:E \to E_1$ over $B$ making the triangle commutate:
 \[
  \xymatrix{& E \ar[d]^g \ar@{.>}[dl]_h \\ E_1 \ar[r]_f & E_2.}
 \]
\end{de}

Formally, we have the following basic properties:

\begin{prop}\label{prop:proj-properties}\
 \begin{enumerate}
  \item Each diffeological vector pseudo-bundle $E_i \to B$ is projective if and only if the direct sum
   $\bigoplus_i E_i \to B$ is projective.
  \item Projectives are closed under taking retracts.
  \item Every bundle subduction to a projective splits smoothly.
 \end{enumerate}
\end{prop}

Recall from~\cite[Section~3.2.5]{CWp} that given a smooth map $X \to B$, we get a diffeological vector 
pseudo-bundle $F_B(X) \to B$.

\begin{lem}\label{lem:free=>proj}
 Let $X \to B$ be a smooth map. The corresponding diffeological vector pseudo-bundle $F_B(X) \to B$ is projective 
 if and only if for every bundle subduction $f:E_1 \to E_2$ over $B$ and any smooth map $g:X \to E_2$ over $B$, 
 there exists a smooth map $h:X \to E_1$ over $B$ such that $g = f \circ h$.
\end{lem}
\begin{proof}
 As usual, this follows from the universal property of $F_B(X) \to B$; see~\cite[Section~3.2.5]{CWp}.
\end{proof}

\begin{prop}\label{prop:plot-proj}
 Every plot $U \to B$ induces a projective diffeological vector pseudo-bundle $F_B(U) \to B$.
\end{prop}
\begin{proof}
 Given any bundle subduction $f:E_1 \to E_2$ over $B$ and any smooth map $g:U \to E_2$ over $B$, 
 we have smooth local liftings $h_i$ of $g$ to $E_1$. Let $\{\lambda_i\}$ be a smooth partition of unity 
 subordinate to the corresponding open cover $\{U_i\}$ of $U$. Then $\sum_i \lambda_i \cdot h_i:U \to E_1$
 is a global smooth lifting of $g$ over $B$, where each $\lambda_i \cdot h_i:U \to E_1$ is defined as 
 \[
  (\lambda_i \cdot h_i)(u) = \begin{cases} \lambda_i(u) h_i(u), & \textrm{if $u \in U_i$} \\
                           \sigma_1 \circ \pi_2 \circ g(u), & \textrm{else} \end{cases}
 \] 
 with $\sigma_1:B \to E_1$ the zero section and $\pi_2:E_2 \to B$ the given diffeological vector pseudo-bundle. 
 The result then follows from Lemma~\ref{lem:free=>proj}.
\end{proof}

As a direct consequence of the above proof, we have:

\begin{cor}\label{cor:global-lift}
 For every bundle subduction, a plot of the total space of the codomain \emph{globally} lifts to a plot of the
 total space of the domain.
\end{cor}

\begin{thm}\label{thm:enough-proj}
 For every diffeological space $B$, the category $\DVPB_B$ has enough projectives, i.e., given any diffeological 
 vector pseudo-bundle $E \to B$, there exists a projective diffeological vector pseudo-bundle $E' \to B$ together 
 with a bundle subduction $E' \to E$ over $B$.
\end{thm}
\begin{proof}
 We take $E' \to B$ to be the direct sum in $\DVPB_B$ of all $F_B(U) \to B$'s indexed over all plots $U \to E$. 
 By Proposition~\ref{prop:plot-proj}, each $F_B(U) \to B$ is projective, and hence by
 Proposition~\ref{prop:proj-properties}(1), 
 $E' \to B$ is projective. By the universal property of $F_B(U) \to B$, we get a bundle map $F_B(U) \to E$ over $B$, 
 and hence a bundle map $E' \to E$ over $B$. By construction, this map is a subduction.
\end{proof}

In summary, for a fixed diffeological space $B$, the pair of projective diffeological vector pseudo-bundles over 
$B$ and the bundle subductions over $B$ forms a projective class.

\subsection{Examples and properties of projectives}\label{subsection:examples}

We first give some examples of projective diffeological vector pseudo-bundles related to 
classical vector bundle theory. To do so, we need:

\begin{lem}\label{lem:preserve-subduction}
 For a smooth map $f:B \to B'$, the pullback $f^*$ sends a bundle subduction over $B'$ to a bundle subduction
 over $B$, and hence it preserves short exact sequences.
\end{lem}
\begin{proof}
 Let $g:E_1' \to E_2'$ be a bundle subduction over $B'$. Then $f^*(E_1') \to f^*(E_2')$ is given by 
 sending $(b,x)$ to $(b,g(x))$. Every plot $p:U \to f^*(E_2')$ gives rise to smooth maps 
 $p_1:U \to B$ and $p_2:U \to E_2'$ via composition with the two projections.
 Since $g$ is a bundle subduction, $p_2$ locally lifts as a smooth map to $E_1'$, which together with 
 $p_1$ induces a local lifting of $p$ to $f^*(E_1')$, showing the first claim.
 
 Since $f^*$ is a right adjoint by Theorem~\ref{thm:adjoint}, it preserves bundle inductions, which together with the first claim
 proves the second one.
\end{proof}

\begin{rem}
 The above lemma also follows from the fact that the pullback $f^*:\DVPB_{B'} \to \DVPB_B$ has a right adjoint $f_{!}$.
 Given a diffeological vector pseudo-bundle $\pi:E \to B$, the bundle $f_{!}(E) \to B'$ is constructed as 
 \[
  f_{!}(E) = \coprod_{b' \in B'} \Gamma(\pi|_{f^{-1}(b')}).
 \]
 When $f^{-1}(b') = \emptyset$, $\Gamma(\pi|_{f^{-1}(b')})$ is $\R^0$.
 A map $p:U \to f_{!}(E)$ is a plot if 
 \begin{enumerate}
  \item the composite $\xymatrix{U \ar[r]^-p & f_{!}(E) \ar[r]^-{\tilde{\pi}} & B'}$ is a plot of $B'$,
   where $\tilde{\pi}$ sends $\Gamma(\pi|_{f^{-1}(b')})$ to $b'$.
  \item for any smooth map $g:V \to U$ and any plot $h:V \to B$ such that the following diagram commute
  \[
   \xymatrix{V \ar[r]^g \ar[d]_h & U \ar[r]^-p & f_{!}(E) \ar[d]^{\tilde{\pi}} \\ B \ar[rr]_f && B',}
  \]
  the map $V \to E$ defined by $v \mapsto (p(g(v)))(h(v))$ is a plot of $E$.
 \end{enumerate}
 It is straightforward to check that $\tilde{\pi}$ is a smooth map between diffeological spaces such that 
 each fibre is a vector space. After dvsification, we get the desired diffeology on the total space $f_{!}(E)$.
 One can check that $f_{!}$ is a functor which is right adjoint to the pullback $f^*$.
 Moreover, each fibre of $f_{!}(E) \to B'$ has the diffeology of the section space; see~\cite[Section~3.1]{CWp}.
 (I would like to thank J. Daniel Christensen for the suggestion of the set-theoretical construction of $f_!(E)$
 in this remark from type theory point of view.)
\end{rem}

Projectiveness is local in the following sense:

\begin{prop}\label{prop:local}
 Let $\pi:E \to B$ be a diffeological vector pseudo-bundle. Assume that there exists a $D$-open cover $\{B_j\}$ of $B$
 such that $i_j^*(E) \to B_j$ is projective in $\DVPB_{B_j}$ for each $j$, where $i_j:B_j \to B$ denotes the inclusion, 
 together with a smooth partition of unity $\{\lambda_j:B \to \R\}$ subordinate to this cover.
 Then $\pi$ is projective in $\DVPB_B$.
\end{prop}
\begin{proof}
 For any bundle subduction $f:E_1 \to E_2$ over $B$ and any bundle map $g:E \to E_2$ over $B$, 
 we get a diagram over $B_j$ for each $j$:
 \[
  \xymatrix{& i_j^*(E) \ar[d]^{i_j^*(g)} \\ i_j^*(E_1) \ar[r]_{i_j^*(f)} & i_j^*(E_2).}
 \]
 Lemma~\ref{lem:preserve-subduction} shows that the horizontal arrow is a bundle subduction over $B_j$.
 By assumption, we have a smooth lifting $h_j:i_j^*(E) \to i_j^*(E_1)$ over $B_j$. 
 Then $\sum_j \lambda_j \cdot h_j:E \to E_1$ is a bundle map over $B$ as we desired.
\end{proof}

We also have the following expected result:

\begin{prop}
 Let $V$ be a projective diffeological vector space, and let $B$ be a smooth manifold. 
 Then the trivial bundle $B \times V \to B$ is projective.
\end{prop}
Surprisingly, note that the result can fail if $B$ is an arbitrary diffeological space; see Example~\ref{ex:not-proj}.
\begin{proof}
 We first reduce the above statement to a special case.
 By Proposition~\ref{prop:local}, it is enough to prove this for the case when $B$ is an open subset
 of a Euclidean space. Recall that every projective diffeological vector space is a retract of direct sums of
 $F(U)$'s for open subsets $U$ of Euclidean spaces (\cite[Corollary~6.15]{W}).
 By Proposition~\ref{prop:proj-properties}(1) and (2), it is enough to show this for the case when $V = F(U)$
 for an open subset $U$ of a Euclidean space.

 Now we prove the statement for the special case when $V = F(U)$, and $B,U$ are Euclidean open subsets.
 As diffeological vector pseudo-bundles over $B$, we have isomorphisms $F_B(B \times U) \cong B \times F(U)$
 of total spaces. The result then follows directly from Proposition~\ref{prop:plot-proj}.
\end{proof}

Combining the above two propositions together with the fact that every fine diffeological vector space 
is projective, we get:

\begin{cor}\label{cor:classical}
 Vector bundles in classical differential geometry are projective.
\end{cor}

However, a projective diffeological vector pseudo-bundle does not need to be locally trivial, 
even when the base space is Euclidean:

\begin{ex}
 Let $f:\R \to \R$ be the square function $x \mapsto x^2$. By Proposition~\ref{prop:plot-proj}, 
 $F_{\R}(\R) \to \R$ is projective. Clearly, the fibre is $\R^0$ for $b < 0$, $\R$ for $b = 0$ 
 and $\R^2$ for $b > 0$. Therefore, a projective diffeological vector pseudo-bundle does not 
 need to be locally trivial.
\end{ex}

Now we discuss some properties of projective diffeological vector pseudo-bundles.

\begin{prop}\label{prop:proj-structure}
 Every projective diffeological vector pseudo-bundle $E \to B$ is a retract of direct sum in $\DVPB_B$ of
 $F_B(U) \to B$ induced by some plots $U \to B$.
\end{prop}
\begin{proof}
 By the proof of Theorem~\ref{thm:enough-proj}, we get a bundle subduction $E' \to E$ over $B$ with $E'$ 
 a direct sum in $\DVPB_B$ of $F_B(U) \to B$ induced by the plots $U \to E$ (and hence some plots $U \to B$, and repetition is allowed).
 Since $E \to B$ is projective, 
 the result then follows from Proposition~\ref{prop:proj-properties}(3).
\end{proof}

Using notations from~\cite{CW19}, we have

\begin{cor}
 Let $E \to B$ be a projective diffeological vector pseudo-bundle. Then $E_b \in \cSV$ for every $b \in B$, 
 i.e., the smooth linear functionals on $E_b$ separate points.
\end{cor}
\begin{proof}
 By Proposition~\ref{prop:proj-structure}, we know that $E$ is a retract of direct sum in $\DVPB_B$ of
 $F_B(U) \to B$ 
 induced by some plots $U \to B$. As $\cSV$ is closed under taking retracts and direct sums
 (\cite[Proposition~3.11]{CW19}), 
 it is enough to show the claim for the special case $F_B(U) \to B$ induced by a plot $p:U \to B$. In this case, 
 the fibre at $b \in B$ is the free diffeological vector space generated by $p^{-1}(b)$
 (\cite[Section~3.2.5]{CWp}), which is a subset of a 
 Euclidean space, and hence $p^{-1}(b) \in \cSD'$, i.e., the smooth functions on $p^{-1}(b)$ separate points.
 The result then follows from~\cite[Proposition~3.13]{CW19}. 
\end{proof}

One would expect that each fibre of a projective diffeological vector pseudo-bundle is a projective 
diffeological vector space. This is equivalent to the statement that the free diffeological vector space 
generated by \emph{any} subset with the subset diffeology of a Euclidean space is projective, 
by a similar argument as above. But I don't know 
whether this is true or not. Nevertheless, we have:

\begin{prop}\label{prop:projective-fibre}
 Let $B$ be a diffeological space. Then every fibre of a projective diffeological vector pseudo-bundle 
 $E \to B$ is a projective diffeological vector space if and only if for every plot $p:U \to B$ and every 
 $b \in B$, the free diffeological vector space generated by $p^{-1}(b)$ is projective.
\end{prop}
\begin{proof}
 ($\Rightarrow$) This follows directly from Proposition~\ref{prop:plot-proj}.

 ($\Leftarrow$) The proof follows from a similar argument as the one in the proof of the above corollary.
\end{proof}

\begin{prop}
 Let $B$ be a discrete diffeological space, i.e., every plot is locally constant.
 Then a diffeological vector pseudo-bundle over $B$ is projective
 if and only if each fibre is a projective diffeological vector space.
\end{prop}
\begin{proof}
 ($\Rightarrow$) This follows from the definition of a discrete diffeological space, together with
 Proposition~\ref{prop:projective-fibre} and~\cite[Corollary~6.4]{W}.

 ($\Leftarrow$) This follows from the fact that every diffeological vector pseudo-bundle over a 
 discrete diffeological space is a coproduct in $\DVPB$ of diffeological vector spaces over a point.
\end{proof}

Also, we have the following results:

\begin{prop}\label{prop:preservation-ses}
 Let $\pi:E \to B$ be a projective diffeological vector pseudo-bundle, and let $\pi_1 \to \pi_2 \to \pi_3$
 be a short exact sequence in $\DVPB_B$, with $\pi_i:E_i \to B$. Then $\Hom_B(\pi,\pi_1) \to \Hom_B(\pi,\pi_2)
 \to \Hom_B(\pi,\pi_3)$ is also a short exact sequence in $\DVPB_B$.
\end{prop}
\begin{proof}
 By Proposition~\ref{prop:proj-structure}, we know that $\pi$ is a retract of direct sum of $F_B(U) \to B$'s
 indexed by some plots $U \to B$. It is straightforward to check that retract and direct product preserve
 short exact sequences in $\DVPB_B$. For the direct product case, one needs Corollary~\ref{cor:global-lift}
 for the subduction part. 
 By the universal property of free bundle induced by a smooth map (\cite[Section~3.2.5]{CWp}), one has a 
 bundle isomorphism over $B$ from
 $\Hom_B(F_B(U),E_i)$ to the set $\Hom_B(U,E_i)$ of all smooth maps $U \to E_i$ preserving $B$, equipped with
 the subset diffeology of $C^\infty(U,E_i)$. Again by Corollary~\ref{cor:global-lift}, it is direct to check
 that $\Hom_B(U,?)$ preserves short exact sequences in $\DVPB_B$. The result then follows by the above
 observations together with the first isomorphism in~\cite[Proposition~3.13]{CWp} 
\end{proof}

\begin{rem}
 The converse of Proposition~\ref{prop:preservation-ses} is false. This is because $\Hom_B(\pi,?)$ always 
 preserves short exact sequences in $\DVPB_B$ for the trivial bundle $\pi:B \times \R \to B$, as it is naturally 
 isomorphic to the identity functor. But the trivial bundle may not be projective; see Example~\ref{ex:not-proj}.
\end{rem}

As a consequence of Proposition~\ref{prop:preservation-ses} and \cite[Proposition~3.12]{CWp}, we have:

\begin{cor}
 If $E_1 \to B$ and $E_2 \to B$ are projective diffeological vector pseudo-bundles, then so is their tensor 
 product $E_1 \otimes E_2 \to B$.
\end{cor}

Since $\Exterior^k E$ is a smooth direct summand of $E^{\otimes k}$ (as a result of~\cite[Lemma~2.11]{P}
and Theorem~\ref{thm:splitses}), by the above corollary and Proposition~\ref{prop:proj-properties}(2), we have:

\begin{cor}
 If $E \to B$ is a projective diffeological vector pseudo-bundle, then so is each exterior product
 $\Exterior^k E \to B$ for $k \geq 1$.
\end{cor}

\subsection{Base change}

\begin{thm}\label{thm:pushforward-proj}
 The pushforward $f_*:\DVPB_B \to \DVPB_{B'}$ sends projectives in the domain to the projectives in the codomain.
\end{thm}
\begin{proof}
 By the adjunction of Theorem~\ref{thm:adjoint}, the following lifting problems are equivalent:
 \[
  \vcenter{\xymatrix{& f_*(E) \ar[d] \ar@{.>}[dl] \\ E_1' \ar[r] & E_2'}}
  \quad \Longleftrightarrow \quad
  \vcenter{\xymatrix{& E \ar[d] \ar@{.>}[dl] \\ f^*(E_1') \ar[r] & f^*(E_2'),}}
 \]
 where $E_1' \to E_2'$ is a bundle subduction over $B'$. By Lemma~\ref{lem:preserve-subduction} and 
 Definition~\ref{de:projective}, we know that the lifting problem on the right has a solution, 
 and hence so is the one on the left.  
\end{proof}

This theorem has several applications. We first give another class of examples of projective 
diffeological vector pseudo-bundles from tangent bundles of diffeological spaces.
To do so, we need the following result:

Note that projective diffeological vector pseudo-bundles are defined in $\DVPB_B$, but they have a similar 
property in $\DVPB$ as follows:

\begin{prop}\label{prop:proj}
 Given a bundle subduction $f:E_1' \to E_2'$ over $B'$ and a bundle map 
 \[
  \xymatrix{E \ar[r]^g \ar[d]_\pi & E_2' \ar[d] \\ B \ar[r]_l & B'}
 \]
 with $\pi$ projective, there exists a bundle map $h:E \to E_1'$ such that $g = f \circ h$.
\end{prop}
\begin{proof}
 By the universal property of pushforward (Proposition~\ref{prop:universal}), we can write $g$ as a bundle map 
 $\tilde{g}:l_*(E) \to E_2'$ over $B'$ followed by the bundle map $\alpha_l:E \to l_*(E)$. 
 By Theorem~\ref{thm:pushforward-proj}, the assumption that $\pi$ is projective over $B$ implies that 
 $\pi_l:l_*(E) \to B'$ is projective over $B'$. Therefore, we have a bundle map $\tilde{h}:l_*(E) \to E_1'$
 over $B'$ such that $\tilde{g} = f \circ \tilde{h}$. Then the composite $\tilde{h} \circ \alpha_l$ is 
 the bundle map $h$ we are looking for.
\end{proof}

Recall from~\cite[Theorem~4.17]{CW16} that every tangent bundle $T^{\dvs}B \to B$ of a diffeological space $B$ 
is a colimit in $\DVPB$ of the tangent bundles $TU \to U$ indexed by the plots $U \to B$. 
Each $TU \to U$ is projective by Corollary~\ref{cor:classical}. 
It is possible that some tangent bundles are projective. (But this is not always the case; 
see Example~\ref{ex:not-proj}.) We show this by an example:

\begin{ex}
 Write $B$ for the cross with the gluing diffeology. We show below that the tangent bundle $T^{\dvs}B \to B$
 is projective. 

 Note that $B$ is the pushout of 
 \[
  \xymatrix{\R & \R^0 \ar[l]_0 \ar[r]^0 & \R}
 \]
 in $\Diff$. It is straightforward to check that the tangent bundle $T^{\dvs} B \to B$ is the colimit of
 \[
  \xymatrix{T \R \ar[d] & T \R^0 \ar[l]_{T0} \ar[r]^{T0} \ar[d] & T \R \ar[d] \\ \R & \R^0 \ar[l]^0 \ar[r]_0 & \R}
 \]
 in $\DVPB$. Write $Tx:T \R \to T^{\dvs} B$ and $Ty:T \R \to T^{\dvs} B$ for the two structural maps.
 Given a bundle subduction $f:E_1 \to E_2$ over $B$ and a bundle map $g:T^{\dvs} B \to E_2$, 
 since $T \R \to \R$ is projective, by Proposition~\ref{prop:proj} we have bundle maps $hx,hy:T \R \to E_1$
 such that $g \circ Tx = f \circ hx$ and $g \circ Ty = f \circ hy$. By the universal property of pushout, 
 we get a desired bundle map $h:T^{\dvs} B \to E_1$ over $B$ with the required property.
\end{ex}

As another consequence of Theorem~\ref{thm:pushforward-proj}, we have the following example which gives 
counterexamples to several arguments:

\begin{ex}\label{ex:not-proj}
 If the free diffeological vector space $F(B)$ is not projective, then the trivial bundle $B \times \R \to B$ 
 is not projective. This happens when the $D$-topology on $B$ is not Hausdorff (\cite[Corollary~3.17]{CW19}).
 The proof of the statement follows from Proposition~\ref{prop:free-geom} and Theorem~\ref{thm:pushforward-proj}.

 This example shows that not every trivial bundle is projective, even when the fibre is a projective (or fine)
 diffeological vector space. It also shows that the pullback functor does \emph{not} preserve projectives, 
 since the trivial bundle $B \times \R \to B$ is the pullback of $\R \to \R^0$ along the map $B \to \R^0$. 
 Furthermore, it shows that not every tangent bundle is projective. For example, $TB \to B$ is not projective 
 when $B$ is an irrational torus, since in this case $TB = B \times \R$ (\cite[combining Examples~3.23 and~4.19(3),
 and Theorem~4.15]{CW16}) and the $D$-topology on $B$ is not Hausdorff.
\end{ex}

Moreover, via Theorem~\ref{thm:pushforward-proj} and Section~\ref{subsection:examples}, 
we get many examples of projective diffeological vector 
spaces from classical differential geometry!

\section{Application to smooth splittings of projective diffeological vector spaces}\label{s:application}

By~\cite[Proposition~3.14 and Theorem~4.2]{CW19}, we know that every finite-dimensional linear subspace 
of a projective diffeological vector space is a smooth direct summand; or in other words, the only 
indecomposable projective diffeological vector space is $\R$. In this section, we use classical smooth
bundle theory and the theory established so far to get some general criteria and interesting examples of 
smooth splittings of projective diffeological vector spaces.

To simplify notation,
we write $V_\pi$ (or $V_E$ when the bundle is understood) for the diffeological vector space obtained from the pushforward of the 
diffeological vector pseudo-bundle $\pi:E \to B$ along the map $B \to \R^0$.

\subsection{General theory}

Here is the general setup.
Given a classical fibre (resp. principal) bundle $E \to B$, we get a linear subduction 
$F(E) \to F(B)$ of diffeological vector spaces which splits smoothly since $F(B)$ is projective. 
We aim to give a bundle-theoretical explanation of its kernel.
In fact, we will prove more general results as follows:

Given a bundle map 
\[
 \xymatrix{E_1 \ar[r]^g \ar[d]_{\pi_1} & E_2 \ar[d]^{\pi_2} \\ B_1 \ar[r]_f & B_2}
\]
from a diffeological vector pseudo-bundle $\pi_1$ to another $\pi_2$, 
by Proposition~\ref{prop:universal}, we get a bundle map $h:f_*(E_1) \to E_2$ over $B_2$
so that $g = h \circ \alpha_f$, where $\alpha_f:E_1 \to f_*(E_1)$ is the structural map introduced at the beginning of Section~\ref{s:pushforward}.
Write $\pi:E \to B_2$ for the kernel of $h$.

Here is the key result:

\begin{thm}\label{thm:kernel}
 Let $(g,f):\pi_1 \to \pi_2$ be a bundle map as above, with $E_1$ locally Euclidean, and $B_2$ Hausdorff and filtered. 
 Then we have a smooth linear map $g_*: V_{\pi_1} \to V_{\pi_2}$ between diffeological vector spaces,
 whose kernel is isomorphic to $V_\pi$ with $\pi:E \to B_2$ defined above.
\end{thm}
\begin{proof}
 By Proposition~\ref{prop:universal}, we get a smooth linear map $g_*: V_{\pi_1} \to V_{\pi_2}$.
 Write $K$ for its kernel. It consists of elements of finite sum $\sum_i e_i$ in $V_{\pi_1}$ with 
 $e_i \in E_1$ such that
 for each $b_2 \in B_2$, the subsum $\sum_{i:\pi_2 \circ g(e_i) = b_2} g(e_i) = 0$.
 So there is a canonical isomorphism $\alpha:V_\pi \to K$ as vector spaces, which is smooth 
 by Lemma~\ref{lem:diffeology-pushforward}.

 Now we use all the extra assumptions to show that the inverse map $\alpha^{-1}$ is smooth. 
 Take a plot $p:U \to K$ and fix $u_0 \in U$. Since the composite $U \to K \hookrightarrow V_{\pi_1}$ 
 is smooth, by Lemma~\ref{lem:diffeology-pushforward}, there exist finitely many plots $p_i:U \to E_1$ by shrinking $U$ around $u_0$
 if necessary, such that $p(u) = \sum_i p_i(u)$ which satisfies that for each $b_2 \in B_2$, the subsum
 $\sum_{i:f \circ \pi_1 \circ p_i(u) = b_2} g(p_i(u)) = 0$ for every $u \in U$. 
 Fix $b_2^0 \in B_2$. Since $B_2$ is Hausdorff, we may assume that the image of the composites
 $f \circ \pi_1 \circ p_i$ do not intersect if their value at $u_0$ are distinct.
 Now take all the index $i$ so that $f \circ \pi_1 \circ p_i(u_0) = b_2^0$, and denote this index subset by $I_{u_0,b_2^0}$.
 Since $E_1$ is locally Euclidean 
 and $B_2$ is filtered, there exist a pointed plot $q:(V,0) \to (B_2,b_2^0)$ and smooth pointed germs
 $h_i:(E_1,p_i(u_0)) \to (V,0)$, so that $q \circ h_i = f \circ \pi_1$
 and $h_i \circ p_i$ is independent of $i$, for all $i \in I_{u_0,b_2^0}$. This then implies that
 $f \circ \pi_1 \circ p_i = q \circ h_i \circ p_i$ are independent of $i$ for all $i \in I_{u_0,b_2^0}$, 
 and hence follows the smoothness of $\alpha^{-1}$. 
\end{proof}

\begin{prop}\label{prop:dvs-subduction}
 If $(g,f):\pi_1 \to \pi_2$ is a bundle subduction, 
 then we get a linear subduction $g_*:V_{\pi_1} \to V_{\pi_2}$ of diffeological vector spaces. 
\end{prop}
\begin{proof}
 This follows directly from Proposition~\ref{prop:universal} and Lemma~\ref{lem:diffeology-pushforward}.
\end{proof}

As a consequence of the above results, we have:

\begin{cor}
 Let $(g,f):\pi_1 \to \pi_2$ be a bundle subduction so that $E_1$ is locally Euclidean, and $B_2$
 is Hausdorff and filtered. Then we have a short exact sequence of diffeological vector spaces
 \[
  0 \to V_{\pi} \to V_{\pi_1} \to V_{\pi_2} \to 0.
 \]
\end{cor}

Now we discuss a special case
\begin{equation}\label{eq:trivial}
 \xymatrix{Y \times \R \ar[d]_{\Pr_1} \ar[r]^{f \times 1_{\R}} & B \times \R \ar[d]^{\Pr_1} \\ Y \ar[r]_f & B,}
\end{equation}
where $f$ is an arbitrary smooth map. 

Observe that

\begin{prop}\label{prop:free}
 The pushforward of $\Pr_1:Y \times \R \to Y$ along $f:Y \to B$ is exactly the free bundle $F_B(Y) \to B$.
\end{prop}
\begin{proof}
 This follows directly from the definition of the free bundle (\cite[Section~3.2.5]{CWp}) and the definition of 
 pushforward of a diffeological vector pseudo-bundle of Section~\ref{s:pushforward}.
\end{proof}

Note that the bundle map $F_B(Y) \to B \times \R$ over $B$ is given by $\sum_i r_i [y_i] \mapsto (b,\sum_i r_i)$, 
where $f(y_i) = b$ for all $i$. We write $\bar{f}_*:\bar{F}_B(Y) \to B$ for its kernel.

\begin{rem}\
 \begin{enumerate}
  \item This proposition generalizes Proposition~\ref{prop:free-geom} by taking $B = \R^0$.

  \item  From above, we know that $F(Y)$ always has a smooth direct summand $\R$ (i.e., $F(Y) \cong \R \oplus \bar{F}(Y)$),
   since $\R$ is a projective diffeological 
   vector space. This can be viewed as a property of the free diffeological vector space, and not every diffeological 
   vector space is free over some diffeological space.

   On the contrary, not every trivial line bundle $B \times \R \to B$ is projective when $B \neq \R^0$ (Example~\ref{ex:not-proj}), so the free 
   bundle $F_B(Y) \to B$ may not have a smooth direct summand $B \times \R \to B$.
 \end{enumerate}
\end{rem}

In the current special case, we have

\begin{cor}\label{cor:kernel}
 Let $f:Y \to B$ be a smooth map, with $Y$ locally Euclidean, and $B$ Hausdorff and filtered.
 \begin{enumerate}
  \item The kernel of $f_*:F(Y) \to F(B)$ is isomorphic to $V_{\bar{f}_*}$ with $\bar{f}_*:\bar{F}_B(Y) \to B$ defined above.
  \item If $f$ is a subduction, then we get a short exact sequence of diffeological vector spaces
   \[
    0 \to V_{\bar{f}_*} \to F(Y) \to F(B) \to 0.
   \]
  \item The pushforward of the free bundle $F_B(Y) \to B$ along $B \to \R^0$ is isomorphic to the free diffeological vector space $F(Y)$.
 \end{enumerate} 
\end{cor}

\begin{rem}
 To make $f_*:F(Y) \to F(B)$ a linear subduction, it is not necessary to require $f:Y \to B$ to be a subduction;
 see Example~\ref{ex:cross}.
\end{rem}

Now we discuss a more special case, which occurs often in practice:
In the diagram~\eqref{eq:trivial}, we further assume that $f$ is a principal $G$-bundle for some diffeological 
group $G$. We give an alternative description of the bundle $V_{\bar{f}_*}$ as follows.

As a setup, assume that $G$ acts smoothly on $Y$ on the right. 
Note that $G$ acts smoothly on $F(G)$ on the left by $G \times F(G) \to F(G)$
given by $g \cdot \sum_i r_i [g_i] = \sum_i r_i [gg_i]$, and it passes to a smooth left action of $G$ on $\bar{F}(G)$,
where $\bar{F}(G)$ is the linear subspace of $F(G)$ consisting of elements of finite sum $\sum_i r_i [g_i]$ with $\sum_i r_i = 0$. 
So we get a commutative square in $\Diff$
\begin{equation}\label{eq:bundlemap}
 \xymatrix{Y \times \bar{F}(G) \ar[r] \ar[d] & \tilde{E} \ar[d]^{\tilde{\pi}} \\ Y \ar[r]_f & B,}
\end{equation}
where $\tilde{E}$ is the quotient of $Y \times \bar{F}(G)$ with $(y,v) \sim (y \cdot g,g^{-1} \cdot v)$ for $y \in Y$, 
$g \in G$ and $v \in \bar{F}(G)$, and $\tilde{\pi}[y,v] = f(y)$.

\begin{lem}
 With the above notations, $\tilde{\pi}$ is a vector bundle over $B$ with fibre $\bar{F}(G)$.
\end{lem}
\begin{proof}
 Let $p:U \to B$ be a plot. 
 Since $f:Y \to B$ is a principal $G$-bundle, we may shrink $U$ so that we have a pullback diagram
 \[
  \xymatrix{U \times G \ar[r]^-\phi \ar[d] & Y \ar[d]^f \\ U \ar[r]_p & B.}
 \]
 We are left to show that there is an isomorphism $\alpha:P \to U \times \bar{F}(G)$ as diffeological vector pseudo-bundles
 over $U$, where $P$ is the pullback of $\xymatrix{U \ar[r]^p & B & \tilde{E}. \ar[l]_{\tilde{\pi}}}$
 We define $\alpha(u,[y,v]) = (u,\theta(u,y) \cdot v)$, where 
 $y = \phi(u,e) \cdot \theta(u,y)$ since $f(y) = p(u) = f(\phi(u,e))$, and $e$ is the identity element in the group $G$. 
 It is clear that $\alpha$ is smooth and fibrewise isomorphic as vector spaces.
 And $\alpha^{-1}$ is given by $(u,v) \mapsto (u,[\phi(u,e),v])$, which is obviously smooth.
\end{proof}

It is straightforward to check that the above square~\eqref{eq:bundlemap} is a bundle map.

\begin{prop}
 Recall that the kernel of the bundle map $F_B(Y) \to B \times \R$ over $B$ is denoted by $\bar{f}_*:\bar{F}_B(Y) \to B$.
 It is isomorphic to $\tilde{\pi}:\tilde{E} \to B$ as vector bundles over $B$.
\end{prop}
\begin{proof}
 The isomorphism as vector bundles over $B$ is given by $\tilde{E} \to \bar{F}_B(Y)$ with 
 $[y,\sum_i r_i [g_i]] \mapsto \sum_i r_i [y \cdot g_i]$, and it is easy to check all the required conditions.
\end{proof}

As a consequence of the above results, we have

\begin{cor}
 Let $f:Y \to B$ be a principal $G$-bundle with $Y$ being locally Euclidean, and $B$ being Hausdorff and filtered. 
 Then we have a short exact sequence of diffeological vector spaces
 \[
  0 \to V_{\tilde{\pi}} \to F(Y) \to F(B) \to 0.
 \]
\end{cor}

Note that when $f:Y \to B$ is a classical fibre (resp. principal) bundle, the conditions ($f$ being a subduction, 
$Y$ locally Euclidean, $B$ Hausdorff and filtered) are satisfied.

\begin{prop}\label{prop:bundle}
 Let $\pi:E \to Y$ be a vector bundle of fibre type a diffeological vector space $V$, 
 and let $f:Y \to B$ be a fibre bundle of fibre type a diffeological space $X$. 
 \begin{enumerate}
  \item\label{covering} If $X$ is finite discrete\footnote{When the fibre of a fibre bundle $f:Y \to B$ is discrete, $f$ is also called a covering},
   then the pushforward $f_*(E) \to B$ is a vector bundle with fibre type $F(X) \otimes V$.
  \item\label{lt} Assume that both $\pi$ and $f$ are locally trivial, 
   and there exists a $D$-open covering $\{B_i\}_i$ of $B$ which 
   trivializes $f$ and simultaneously the $D$-open covering $\{f^{-1}(B_i)\}_i$ trivializes $\pi$. Then the pushforward 
   $f_*(E) \to B$ is also a locally trivial vector bundle of fibre type $F(X) \otimes V$.
  \item\label{fb} If $\pi$ is trivial, then $f_*(E) \to B$ is a vector bundle of fibre type $F(X) \otimes V$.
 \end{enumerate}
\end{prop}
\begin{proof}
 \eqref{covering} Let $p:U \to B$ be a plot. Since $f:Y \to B$ is a covering with fibre type $X$, we may shrink
  $U$ to get a pullback diagram
  \[
   \xymatrix{U \times X \ar[r]^-\phi \ar[d] & Y \ar[d]^f \\ U \ar[r]_p & B.}
  \]
  Since $\pi:E \to Y$ is a vector bundle of fibre type $V$, for each $x \in X$, we may further shrink $U$ to get 
  a pullback diagram
  \[
   \xymatrix{U \times \{x\} \times V \ar[r]^-{\psi_x} \ar[d] & E \ar[d]^\pi \\ U \times \{x\} \ar[r]_-{\phi|_{U \times \{x\}}} & Y.}
  \]
  As $X$ is finite discrete, we gather these together and get a pullback diagram
  \[
   \xymatrix{U \times X \times V \ar[r]^-\psi \ar[d] & E \ar[d]^\pi \\ U \times X \ar[r]_-\phi & Y.}
  \]
  Write $P$ for the pullback of $\xymatrix{U \ar[r]^p & B & f_*(E). \ar[l]}$
  Then $P$ consists of elements of the form $(u,\sum_i e_{y_i})$ with $p(u) = f(y_i)$ for all $i$.
  Define $U \times (F(X) \otimes V) \to P$ by linear expansion of $(u,[x] \otimes v) \mapsto (u,\psi(u,x,v))$.
  It is straightforward to check that this map is smooth and an isomorphism of vector spaces, and its inverse is also smooth.

 \eqref{lt} and~\eqref{fb} can be proved in a similar way.
\end{proof}

\begin{cor}
 If $f:Y \to B$ is a (locally trivial) fibre bundle of fibre type a diffeological space $X$, 
 then $F_B(Y) \to B$ is a (locally trivial) vector bundle of fibre type $F(X)$.
\end{cor}
\begin{proof}
 This follows immediately from Propositions~\ref{prop:free} and~\ref{prop:bundle}.
\end{proof}

\subsection{Examples}

Now we deal with the case of principal bundle whose group $G$ is discrete. 
In this case, $F(G)$ is a fine diffeological vector space whose dimension matches the cardinality of $G$, 
and $\bar{F}(G)$ is a codimension-one linear subspace of $F(G)$, and hence also a fine diffeological vector space.

\begin{ex}\label{ex:RPn}
 For the principal $\Z/2\Z$-bundle $S^n \to \R P^n$, $F(\Z/2\Z) \cong \R^2$ and $\bar{F}(\Z/2\Z) \cong \R$.
 And therefore, the bundle $\tilde{\pi}$ in the commutative square~\eqref{eq:bundlemap} in the previous subsection
 can be viewed as the quotient of $S^n \times \R$ with the equivalence relation given by $(z,x) \sim (-z,-x)$, which
 is the tautological line bundle $\gamma^1_n$ on $\R P^n$. So we have an isomorphism
 \begin{equation}\label{eq:SnRPn}
  F(S^n) \cong F(\R P^n) \oplus V_{\gamma^1_n}.
 \end{equation}
 Taking $n=1$, $\gamma^1_1$ is the M\"obius band. Moreover, since $\R P^1$ is diffeomorphic to $S^1$, we get
 \begin{equation}\label{eq:decomposition}
  F(S^1) \cong F(S^1) \oplus V_{\gamma^1_1} \cong \ldots \cong F(S^1) \oplus (V_{\gamma^1_1})^m
 \end{equation}
 for any $m \in \N$.
\end{ex}

By some results from~\cite{MS}, we have

\begin{ex}\
 \begin{enumerate}
  \item Since the tangent bundle $TS^n \to S^n$ direct sum the normal bundle (which is the trivial line bundle) of $S^n$ in $\R^{n+1}$
   is a trivial bundle over $S^n$ of rank $n+1$, we get 
   \[
    F(S^n)^{n+1} \cong F(S^n) \oplus V_{TS^n}.
   \]
   Moreover, by~\cite{A}, $V_{TS^n}$ has a smooth direct summand $F(S^n)^{\rho(n+1) - 1}$, where $\rho(n+1) = 2^c + 8d$
   with $n+1 = 2^b(2a+1)$, $b = c + 4d$ and $0 \leq c \leq 3$.

  \item Since the tangent bundle $T \R P^n \to \R P^n$ direct sum the trivial line bundle over $\R P^n$ is 
   isomorphic to the direct sum of $(n+1)$-copies of the tautological line bundle $\gamma^1_n \to \R P^n$, we get
   \[
    (V_{\gamma^1_n})^{n+1} \cong F(\R P^n) \oplus V_{T \R P^n}.
   \]

  \item The total space of the tangent bundle $TS^n \to S^n$ can be viewed as a submanifold of
   $\R^{n+1} \times \R^{n+1}$, with the first component for the base and the second one for the tangent part.
   If we identify $(x,v)$ with $(-x,-v)$ in $TS^n$, we get the total space of the 
   tangent bundle $T \R P^n \to \R P^n$; if we identify $(x,v)$ with $(-x,v)$ in $TS^n$, we get another locally 
   trivial vector bundle $\pi:E \to \R P^n$ of rank $n$. (In the case $n=1$, $\pi$ is exactly the M\"obius band
   over $\R P^1$; notice the difference from Example~\ref{ex:RPn}, based on the different meaning of the coordinates!) Write $f:S^n \to \R P^n$ for the quotient map. 
   Note that $E \to f_*(TS^n)$ given by $[x,v] \mapsto (x,v)+(-x,v)$ is a bundle map over $\R P^n$, using 
   Proposition~\ref{prop:bundle}\eqref{covering}, which is the kernel of the canonical bundle map $f_*(TS^n) \to T \R P^n$. Hence, we have an isomorphism
   \[
    V_{TS^n} \cong V_{T \R P^n} \oplus V_\pi,
   \]
   which also recovers the first isomorphism in~\eqref{eq:decomposition} in Example~\ref{ex:RPn}.
 \end{enumerate}

 Therefore, if we combine the three isomorphisms in this example, we get
 \[
  F(\R P^n) \oplus F(S^n)^{n+1} \cong F(S^n) \oplus V_\pi \oplus (V_{\gamma^1_n})^{n+1}.
 \]
 By taking $n=1$, we obtain
 \[
  F(S^1)^3 \cong F(S^1) \oplus (V_{\gamma^1_1})^3.
 \]
\end{ex}

Finally, we show by the following example that the extra condition of filteredness added to the results in the previous subsection 
is necessary.

\begin{ex}
 Let $\Z/2 \Z$ act on $\R$ by $\pm 1 \cdot x = \pm x$, and write $B$ for the quotient space.
 Then $B$ is weakly filtered but not filtered (\cite[Example~4.7]{CW17}), and $B$ with the $D$-topology is
 homeomorphic to the subspace $[0,\infty)$ of $\R$ (hence Hausdorff).
 Write $f:\R \to B$ for the quotient map, and write $K$ for the kernel of $F(\R) \to F(B)$.
 It consists of elements of the form of finite sum $\sum_i r_i [x_i]$ with $r_i,x_i \in \R$ such that for 
 every fixed $x \in X$, the subsum $\sum_{i:x_i = \pm x} r_i = 0$. So, $p:\R \to K$ defined by $t \mapsto [t] - [-t]$ is a plot of $K$.
 On the other hand, the map $f_*:F_B(\R) \to B$ has fibre $\R$ over $[0] \in B$ and fibre $\R^2$ over $[b] \in B$
 for $b \neq 0$. Hence, $\bar{f}_*:\bar{F}_B(\R) \to B$ has fibre $\R^0$ over $[0] \in B$ and fibre $\R$ over 
 $[b] \in B$ for $b \neq 0$. The canonical smooth linear bijection $\alpha:V_{\bar{f}_*} \to K$ is not an isomorphism 
 of diffeological vector spaces since $\alpha^{-1} \circ p$ is not a plot of $V_{\bar{f}_*}$. If it were, 
 then by iterated use of Lemma~\ref{lem:diffeology-pushforward} there exist finitely many smooth germs
 $(p^1_{i,j},p^2_{i,j}):\R \to \R_{(\text{base})} \times \R_{(\text{fibre})}$ at $0 \in \R$ such that 
 \[
  p(t) = \alpha( \sum_{i,j} \alpha_g(\alpha_f (p^1_{i,j}(t),p^2_{i,j}(t))) ),
 \]
 where $g:B \to \R^0$, both $\alpha_f$ and $\alpha_g$ are structural maps from Section~\ref{s:pushforward},
 the range of $j$ depends on $i$, $f \circ p^1_{i,j}$ is independent of $j$ for any 
 fixed $i$, $p^2_{i,j}(t) = 0$ whenever $p^1_{i,j}(t) = 0$ (by the description of $V_{\bar{f}_*}$, which causes the contradiction as follows). 
 By evaluating at $t=0$, we know that $\sum_{i,j:p^1_{i,j}(0) = x} p^2_{i,j}(0) = 0$ for any fixed $x \in \R \setminus \{0\}$. By continuity of 
 the $p^2_{i,j}$'s, we know that $\sum_{i,j:p^1_{i,j}(t)=t} p^2_{i,j}(t) \neq 1$ for $t \neq 0$ but sufficiently close to $0$, 
 which implies that $\alpha^{-1} \circ p$ cannot be a plot. 
\end{ex}

\vspace*{10pt}

\end{document}